\documentclass[12pt]{article}
\usepackage{amssymb,amsmath,amsthm,rotating}
\newcommand{\n}{\noindent}
\newcommand{\bb}[1]{\mathbb{#1}}
\newcommand{\cl}[1]{\mathcal{#1}}

\newcommand{\vp}{\varepsilon}

\theoremstyle{plain}

\newtheorem{thm}{Theorem}

\theoremstyle{definition}

\begin{document}

\title{Variations on a theme of Beurling}

\author{Ronald G. Douglas\footnote{During the time the research leading up to this note was carried out,  the author was partially supported by a grant from the National Science Foundation.}}

\date{}
\maketitle

\thispagestyle{empty}

\abstract{
Interpretations of the Beurling--Lax--Halmos Theorem on invariant subspaces of the unilateral shift are explored using the language of Hilbert modules. Extensions and consequences are considered  in  both the one and multivariate cases with an emphasis on the classical Hardy, Bergman  and Drury--Arveson spaces.
}
 
\vfill

\n {\bf 2000 Mathematics Subject Classification}:\ 46E22, 46J10, 47A15.

\n {\bf  Key words and phrases}. Invariant subspaces, Beurling's Theorem, Hilbert modules.

\newpage 

\setcounter{section}{-1}

\section{Introduction}\label{sec0}

\indent 

In a classic paper \cite{B}, Beurling posed and answered two fundamental questions for the unilateral shift operator on Hilbert space and its adjoint. The first problem was the characterization of the cyclic vectors for the forward shift operator, while the second one concerned the spanning of an invariant subspace for the backward shift by its eigenvectors or, more generally,  its generalized eigenvectors. To obtain these solutions, he recast the questions into the language of   function theory and then  recalled results of Nevanlinna  and Riesz on the inner-outer factorization of functions in the Hardy space on the unit disk and the structure of inner functions. In particular, a vector is cyclic for the unilateral shift if and only if it is outer or its inner factor is a scalar and spectral synthesis holds for an invariant subspace for the backward shift if and only if the inner function, representing its orthogonal complement, has no singular inner factor and the zeros of the Blaschke product have multiplicity one. For higher multiplicity zeros one must also include generalized eigenvectors. Thus Beurling's solutions rested on the results in function theory obtained a decade or two earlier.

The result from his paper, which is  the best remembered, is the representation of invariant subspaces of the unilateral shift in terms of inner functions or, what is usually called Beurling's Theorem. Perhaps what is somewhat surprising is that this result is really a statement about the structure of isometries and could have been obtained as a corollary to von~Neumann's result \cite{vN} two decades earlier, what is now usually called the Wold decomposition \cite{W}. This fact becomes transparent if one adopts a Hilbert module point of view. In this note we will do that examining various interpretations and generalizations of Beurling's results in the context of Hilbert modules of holomorphic functions on domains in ${\bb C}^m$ such as the unit ball ${\bb B}^m$ and the polydisk ${\bb D}^m$ in ${\bb C}^m$ for $m\ge 1$. Many of these ideas occurred and were developed by the author in collaboration with Jaydeb Sarkar.

\section{Preliminaries}\label{sec1}

\indent 

We will restrict definitions to the context needed in this note. For a more detailed presentation of Hilbert modules see \cite{DM}, \cite{DP}. 

A Hilbert module ${\cl H}$ over ${\bb C}[\pmb{z}], \pmb{z} = (z_1,\ldots, z_m)$ with $m\ge 1$, is a Hilbert space ${\cl H}$ and a unital module action
\[
 {\bb C}[\pmb{z}]\times {\cl H}\to {\cl H}
\]
such that each operator $M_p$ for $p$ in ${\bb C}[\pmb{z}]$ defined $M_pf= p\cdot f$ for $f$ in ${\cl H}$ is bounded.

Examples of Hilbert modules are the Hardy space $H^2({\bb B}^m)$ on the unit ball ${\bb B}^m$ in ${\bb C}^m$ for $m\ge 1$, the Hardy space $H^2({\bb D}^m)$ on the polydisk ${\bb D}^m$ in ${\bb C}^m$ for $m\ge 1$, the Bergman space $L^2_a({\bb B}^m)$ on ${\bb B}^m$ and also the Bergman space $L^2_a({\bb D}^m)$ on ${\bb D}^m$ for $m\ge 1$ and many more. 

Not all Hilbert modules can be represented in a natural way as Hilbert spaces of holomorphic functions in which module multiplication agrees with the pointwise multiplication of functions, but we will focus in this note on those that can.

Recall that $H^2({\bb B}^m)$ can be identified as the closure of ${\bb C}[\pmb{z}]$ in $L^2(\partial{\bb B}^m)$ for Lebesgue measure on $\partial {\bb B}^m$. Moreover, $H^2({\bb D}^m)$ is the closure of ${\bb C}[\pmb{z}]$ in $L^2((\partial{\bb D})^m)$ for the product measure on $(\partial{\bb D})^m$. Finally, the Bergman spaces are the closures of ${\bb C}[\pmb{z}]$ in $L^2({\bb B}^m)$ and $L^2_a({\bb D}^m)$, respectively,  for Lebesgue measure on ${\bb B}^m$ and ${\bb D}^m$, respectively. Module multiplication in all of these cases is defined by the pointwise multiplication of functions.

If ${\cl H}$ is a Hilbert module over ${\bb C}[\pmb{z}]$, then there is a natural way to make the Hilbert space tensor product, ${\cl H}\otimes {\cl E}$, into a Hilbert module over ${\bb C}[\pmb{z}]$ for each coefficient Hilbert space ${\cl E}$. One defines $p\cdot (f\otimes x) = (M_p\otimes I_{\cl E})(f\otimes x) = (p\cdot f)\otimes x$ for $f$ in ${\cl H}$, $x$ in ${\cl E}$ and $p$ in ${\bb C}[\pmb{z}]$. This construction enables one to increase the multiplicity of a Hilbert module. In particular, $H^2({\bb D})\otimes {\cl E}$ is the Hardy space on ${\bb D}$ with multiplicity equal to $\dim{\cl E}$.

The result of Beurling was generalized by Lax \cite{L} and Halmos \cite{H} to obtain the following theorem here stated in the language of submodules of Hilbert modules.

\begin{thm}[Beurling--Lax--Halmos]\label{thm1}
Let ${\cl S}$ be a non-zero submodule of $H^2({\bb D})\otimes{\cl E}$ for some Hilbert space ${\cl E}$. Then there exists a subspace ${\cl E}_*$ of ${\cl E}$ such that ${\cl S}$ and $H^2({\bb D})\otimes {\cl E}_*$ are unitarily equivalent Hilbert modules. 
\end{thm}

Recall that Hilbert modules ${\cl H}_1$ and ${\cl H}_2$ over ${\bb C}[\pmb{z}]$ are said to be unitarily equivalent if there exists a unitary module map $U\colon \ {\cl H}_1\to {\cl H}_2$; that is, a unitary map $U$ such that $p\cdot Uf = U(p\cdot f)$ for $f$ in ${\cl H}$ and $p$ in ${\bb C}[\pmb{z}]$.

The BLH result follows directly from von~Neumann's result in \cite{vN}.

\begin{thm}[von Neumann]\label{thm2}
Every isometry on Hilbert space is unitarily equivalent to an operator of the form $(M_z \otimes I_{\cl D})\oplus V$ for some Hilbert space ${\cl D}$ and unitary operator $V$.
\end{thm}

Moreover, the proof of the BLH Theorem follows from that of von~Neumann by proving that $\dim {\cl D} \le \dim {\cl E}$ for ${\cl S} \subseteq {\cl H}^2({\bb D})\otimes {\cl E}$ and that  there is no unitary $V$ in the representation. 
Further, one obtains the usual representation for ${\cl S}$ by using the fact that the module map $U\colon \ H^2({\bb D})\otimes {\cl E}_*\to {\cl S} \subseteq  H^2({\bb D})\otimes {\cl E}$, has the form $(Uf)(z) = \Theta(z)f(z)$ for some holomorphic map $\Theta\colon \ {\bb D}\to {\cl L}({\cl E}_*,{\cl E})$. Since $\Theta$ is a contraction, we have $\|\Theta(z)\|\le 1$ for $z$ in ${\bb D}$, which implies that $\Theta$ has radial limits $\Theta(e^{it})$ on ${\bb T}= \partial{\bb D}$ a.e., which are isometric a.e.

Hence the Hardy module on ${\bb D}$ has the property that all non-zero submodules of the higher multiplicity version have the same form. If we take the Hardy space on a domain $\Omega \subset {\bb C}$ with $\partial\Omega$ a simple closed curve, it will have the same property regarding submodules of its higher multiplicity versions (see \cite{AD}). One can ask if  there are any other Hilbert modules with this property?

If the Hilbert module ${\cl H}$ has no proper submodule, then ${\cl H}$ would satisfy this criterion for trivial reasons. However,  here we eliminate such a possibility by   focusing on the case of quasi-free Hilbert modules \cite{DM} which consist of holomorphic functions on some domain in ${\bb C}^m$ for $m\ge 1$; that is, ${\cl H}\subseteq  \text{hol}(\Omega,{\cl E})$ for a bounded domain such as $\Omega = {\bb B}^m$ or ${\bb D}^m$ and a Hilbert space ${\cl E}$. We assume ${\cl H}$ is the closure of the algebraic tensor product ${\bb C}[\pmb{z}]\otimes {\cl E}$. Such a space is a reproducing kernel Hilbert space, where the kernel $K(\pmb{z},\pmb{w})\colon \ \Omega\times \Omega\to {\cl L}({\cl E})$ is defined such that $K(\pmb{z},\pmb{w})  = E_{\pmb{z}} E^*_{\pmb{w}}$, where $E_{\pmb{z}}\colon \ \Omega\to {\cl E}$ is evaluation at $\pmb{z}$ in $\Omega$, which is bounded.  One says that $\dim {\cl E}$ is the multiplicity of ${\cl H}$.

\section{Quasi-Free Hilbert Modules of Multiplicity One}\label{sec2}

\indent 

For the class of quasi-free Hilbert modules  of multiplicity one, one can decide for which modules ${\cl R}$ all submodules of ${\cl R}$ are isometrically isomorphic to ${\cl R}$.

\begin{thm}\label{thm3}
Suppose ${\cl R}$ is a quasi-free Hilbert module over ${\bb B}^m$ of multiplicity one such that each submodule ${\cl S}$ of ${\cl R}$ is isometrically isomorphic to ${\cl R}$. Then $m=1$ and ${\cl R}$ is isometrically isomorphic to $H^2({\bb D})$ and $M_z$ on ${\cl R}$ is the Toeplitz operator $T_\varphi$, where $\varphi$ is a conformal self map of ${\bb D}$ onto itself. 
\end{thm}

\begin{proof}
Suppose ${\cl R}$ is a quasi-free Hilbert module of multiplicity one over ${\bb C}[\pmb{z}]$ such that every proper submodule is isometrically isomorphic to ${\cl R}$. Then note  that ${\cl S} = \{f\in {\cl R}\colon \ f(\pmb{0}) = 0\}$ is a proper submodule of ${\cl R}$ of codimension one. By hypothesis, the results of \cite{DS} apply. Hence $m=1$ and the module isometry between ${\cl R}$ and ${\cl S}$ yields the identification of ${\cl R}$ and $H^2({\bb D})$, and $M_z$ must be a Toeplitz operator on ${\bb D}$. Moreover, identification of the point spectrum of  $T^*_\varphi$ on ${\bb D}$ and the index completes the proof.
\end{proof}

The same result holds if one replaces ${\bb B}^m$ by $(\partial{\bb D})^m$. An earlier result of Richter \cite{Rich} showed that no proper submodule ${\cl S}$ of the Bergman module $L^2_a({\bb D})$ is isometrically isomorphic to $L^2_a({\bb D})$ revealing that this module has the opposite property for submodules. Actually, one can show, as was established in \cite{DS}, that this statement holds for most subnormal Hilbert modules. The proof there depends on the maximum principle.

However, the following result, based on an operator theoretic approach, covers most cases of interest. Recall that a multivariate weighted shift is defined on the Hilbert space $\ell^2({\bb Z}^m_+)$ by a multi-sequence of weights $\Lambda=\{\lambda_{\pmb{\alpha}}\}_{\pmb{\alpha}\in A}$, where $M\ge \lambda_\alpha>0$ for some positive $M$ and $A = ({\bb Z}^m_+)$, such that the coordinate operators are defined $M_ie_{\pmb{\alpha}} = \lambda_{\pmb{\alpha}} e_{\pmb{\alpha} +\delta_i}$ for $e_{\pmb{\alpha}}$ in $\ell^2({\bb Z}^m_+)$ and $\pmb{\alpha} + \delta_i = (\alpha_1,\ldots, \alpha_i+1,\ldots, \alpha_m)$ for $\pmb{\alpha}$ in $A$ and $i=1,2,\ldots, m$. The multi-shift will be said to be \emph{strictly hyponormal} if $\lambda_{\pmb{\alpha}+\delta_i} - \lambda_{\pmb{\alpha}}
>0$ for $\pmb{\alpha}$ in $A$ and $i=1,\ldots, m$.

\begin{thm}\label{thm4}
If $\Lambda$ is a family of weights defining a strictly hyponormal multi-shift and ${\cl S} \subset \ell^2({\bb Z}^m_+)$ is a submodule isometrically isomorphic to $\ell^2({\bb Z}^m_+)$, then ${\cl S} = \ell^2({\bb Z}^m_+)$.
\end{thm}

\begin{proof}
Let $W$ be the unitary module map between $\ell^2({\bb Z}^m_+)$ and ${\cl S}\subseteq \ell^2({\bb Z}^m_+)$. If $f=We_{\pmb{0}}$, then an easy argument applied to the expansion, $f = \sum\limits_{\pmb{\alpha}\in A} a_{\pmb{\alpha}} e_{\pmb{\alpha}}$, and the actions of the coordinate multipliers $M_{z_i}$, for $i=1,\ldots, m$, shows that $a_{\pmb{\alpha}} = 0$ for $\pmb{\alpha} \ne \pmb{0}$.
\end{proof}

A careful but straightforward modification of this argument, given in \cite{DS}, is shown to apply to the Drury--Arveson space $H^2_m$ obtaining the same result on submodules of $H^2_m$.
Recall that $H^2_m$ can be identified as the symmetric Fock space (see \cite{Arv}).

\section{The Question of Multiplicity}\label{sec3}

\indent

The Beurling--Lax--Halmos Theorem allows one to say more about multiplicity in the case of the Hardy module. In particular, suppose $f$ is in $H^2({\bb D})\otimes {\cl E}$ and $[f]$ denotes the submodule of $H^2({\bb D})\otimes {\cl E}$ generated by $f$. Then $[f]\cong H^2({\bb D})$ or every ``multiplicity one'' (singly-generated) submodule looks like $H^2({\bb D})$. This is not true in general even if one relaxes the requirement as the following example shows.

\begin{thm}\label{thm5}
Consider the vector $1\oplus z$ in $L^2_a({\bb D})\oplus L^2_a({\bb D}) \cong L^2_a({\bb D}) \otimes {\bb C}^2$ and let $[1\oplus z]$ denote the cyclic submodule it generates. Then for no submodule ${\cl S}\subseteq L^2_a({\bb D})$ is ${\cl S}\cong [1\oplus z]$.
\end{thm}

\begin{proof}
Suppose  $W$ is a module isomorphism from $[1\oplus z]$ onto ${\cl S} \subseteq L^2_a({\bb D})$ for some submodule ${\cl S}$. If $f=W(1\oplus z)$, then the  facts that the functions are continuous and the closed support of Lebesgue measure on ${\bb D}$ is ${\bb D}$ implies that $|f(z)|^2 = 1+|z|^2$ for $z$ in ${\bb D}$. Using the Taylor series expansion one can show this is impossible for any holomorphic function $f$ on ${\bb D}$.
\end{proof}

There is  considerable literature, going back at least to Polya, on the question of when the absolute value of a polynomial can be represented as the norm of a vector-valued polynomial or vice versa (see \cite{DA}). These results are related to the theorem although the proof of this special case requires only the uniqueness of Taylor coefficients in expansion of a function in terms of $z$ and $\bar z$. It seems likely that some interesting results could emerge from applying this classical theory to the context of Hilbert modules.

The theorem shows we cannot identify cyclic submodules of $L^2_a({\bb D})\otimes {\cl E}$ with submodules of $L^2_a({\bb D})$ even given the great variety of the cyclic submodules of $L^2_a({\bb D})$ that are known to exist (see \cite{ABF}). However, a reinterpretation of a result of Trent and Wick \cite{TW}, given in \cite{DS2}, shows that some aspects of this property persist for the Hardy modules on ${\bb B}^m$ and ${\bb D}^m$.

Let $A({\bb B}^m)$ and $A({\bb D}^m)$, respectively, denote the closure of ${\bb C}[\pmb{z}]$ in the supremum norm on ${\bb B}^m$ or ${\bb D}^m$, respectively.

\begin{thm}\label{thm6}
Suppose $\{\psi_i\}^N_{i=1}$ are vectors in $A({\bb B}^m)$ or $A({\bb D}^m)$, respectively, and $[\psi_1\oplus\cdots\oplus \psi_N]$ is the cyclic submodule in $H^2({\bb B}^m) \otimes {\bb C}^N$ or $H^2({\bb D}^m)\otimes {\bb C}^N$, respectively, that it generates. Then there exists a vector $f$ in $H^2({\bb B}^m)$ or  an $f$ in $H^2({\bb D}^m)$, respectively, such that
\[
 [\psi_1\oplus\cdots\oplus \psi_N] \cong [f].
\]
In particular, $[\psi_1 \oplus\cdots\oplus \psi_n]$ is isometrically isomorphic to a submodule of $H^2({\bb B}^m)$ or $H^2({\bb D}^m)$, respectively.	
\end{thm}

\begin{proof}
The question comes down to the existence of a holomorphic function $f$ on ${\bb B}^m$ or ${\bb D}^m$, respectively, such that $|f(\pmb{z})|^2 = \sum\limits^n_{i=1} |\psi_i(\pmb{z})|^2$ for $\pmb{z}$ in $\partial{\bb B}^m$ a.e.\ or in $(\partial {\bb D})^m$ a.e. This is a classical problem in the function theory of several complex variables with an affirmative answer in this case (see \cite{Ru}).
\end{proof}

This result holds more generally for $\{\psi_i\}^N_{i=1} \subset H^2({\bb B}^m)$ or in H$^2({\bb D}^m)$, respectively, so long as all of the quotients, $\psi_i(z)/\psi_j(z)$, $1\le i,j\le N$, are continuous on $\partial{\bb B}^m$ or on $(\partial{\bb D})^m$, respectively. These arguments can be turned around to show the equivalence of the module isomorphism and the representation of the absolute value of the functions on the boundary. However, the result in function theory on which these results are based (see \cite{TW}) are thought to be false for general functions in $H^2({\bb B}^m)$, or $H^2({\bb D}^m)$, respectively (see \cite{Ru}). The latter, if correct, would seem strange since it means the answer to this operator theoretic question of ``multiplicity'' rests on the relative boundary behavior of the functions.

We formalize these ideas as follows.
\medskip 

\n {\bf Question 1.} For $\{\psi_i\}^N_{i=1}$ in $H^2({\bb B}^m)$ or $H^2({\bb D}^m)$, respectively, does there exist $f$ in $H^2({\bb B}^m)$ or $H^2({\bb D}^m)$, respectively, such that $|f(z)|^2 = \sum\limits^N_{i=1}|\psi_i(z)|$ for $z$ in $\partial {\bb B}^m$ a.e.\ or $z$ in $(\partial{\bb D})^m$ a.e.  Equivalently, does there exist a submodule ${\cl L}$ of $H^2_a({\bb B}^m)$ or of $H^2({\bb D}^m)$, respectively, such that ${\cl L}\cong [\psi_1\oplus\cdots\oplus \psi_N]$?

For the Bergman spaces, we formulate a somewhat related question.\medskip 

\n {\bf Question 2.} For $\{\psi_i\}^N_{i=1}$ in $L^2_a({\bb B}^m)$ or $L^2_a({\bb D}^m)$, respectively, what is the smallest $k\ge 1$ such that the Hilbert module $[\psi_1 \oplus \cdots \oplus \psi_N]$ is isometrically isomorphic to a submodule of $L^2_a({\bb B}^m)\otimes {\bb C}^k$ or $L^2_a({\bb D}^m) \otimes {\bb C}^k$, respectively?

Although both questions are framed in the language of Hilbert modules, they are equivalent to questions concerning the absolute values of holomorphic functions in several complex variables.

\section{Cyclic Submodules}\label{sec4}

\indent

Note that the earlier theorem shows that the cyclic submodules of $H^2({\bb B}^m)$ or $H^2({\bb D}^m)$, respectively, are not all isomorphic. More is true for cyclic submodules of $L^2_a({\bb B}^m)$ and $L^2_a({\bb D}^m)$.

\begin{thm}\label{thm7}
Let $f_1$ and $f_2$ be vectors in $L^2_a({\bb B}^m)$ or $L^2_a({\bb D}^m)$, respectively, for $m\ge 1$ such that $[f_1]  \cong [f_2]$. Then ${\cl Z}(f_1) = {\cl Z}(f_2)$, where ${\cl Z}(f_i) = \{\pmb{z}\in {\bb B}^m\colon \ f_i(\pmb{z}) = 0\}$ or $\{\pmb{z}\in {\bb D}^m\colon \ f_i(\pmb{z}) = 0\}$, respectively.
\end{thm}

\begin{proof}
Using the fact that the polynomials spanned by the monomials $z^\alpha\bar z^\beta$ for $\pmb{\alpha},\pmb{\beta}$ in $A$ are dense in $C(\text{clos }  {\bb B}^m)$, one can show that $|f_1|^2 dV = |f_2|^2dV$ a.e., where $V$ is volume measure on ${\bb B}^m$ or ${\bb D}^m$, respectively (see \cite{DS}). Thus, one has ${\cl Z}(f_1) = {\cl Z}(f_2)$ since $f_1$ and $f_2$ are continuous on ${\bb B}^m$ or ${\bb D}^m$, respectively, and the closed support of the volume measure is the closed domain. This completes the proof.
\end{proof}

This result provides an uncountable family of nonisometrically isomorphic cyclic submodules of $L^2_a({\bb B}^m)$ or $L^2_a({\bb D}^m)$, respectively, for $m>1$ by choosing a family of functions with different zero sets. For example, consider the  family $f_{\pmb{a}}(\pmb{z}) = \pmb{a}= \sum^m_{i=1} a_iz_i$ for $\pmb{a}$ in ${\bb C}^m$ and $\|\pmb{a}\|=1$.

One can say more about the zero varieties  ${\cl Z}(f_1)$ and ${\cl Z}(f_2)$ but we won't pursue that here.
The question of a converse concerns the nature of the quotients of $f_1$ and $f_2$ when the zero varieties are equal, which as one knows, can be quite complicated.

The above result extends to other subnormal modules on ${\bb B}^m$, so long as the closed support of the measure equals $\text{clos}({\bb B}^m)$. (Actually one can do with much less such as the  space defined  as the closure of ${\bb C}[\pmb{z}]$ in the $L^2$-space for volume measure on $\{\pmb{z}\in {\bb B}^m\colon \ \|\pmb{z}\| > \varepsilon\}$ for $0<\vp < 1$.)

Analogous questions to those discussed above but for the Hardy space $H^2({\bb B}^m)$ have very different answers. For example, if $f$ is an inner function on ${\bb D} = {\bb B}^1$, then $[f] \cong [1] = H^2({\bb D})$. But ${\cl Z}(f)$ is nonempty unless $f$ is a singular inner function and ${\cl Z}(1)=\emptyset$. Similarly, because nontrivial inner functions exist on ${\bb B}^m$, we have the same phenomenon there.

\section{Similarity of Cyclic Submodules}\label{sec5}

\indent 

The Rigidity Theorem in \cite{DP} shows for ideals $I$ in ${\bb C}[\pmb{z}]$ satisfying certain properties, the closure of two ideals $I_1$ and $I_2$ in  a quasi-free Hilbert module are similar if and only if the ideals coincide. Unless $m=1$, principal ideals do not satisfy the additional assumptions and hence there is no general similarity result in this case. 

Let us raise a question about similarity in the simplest possible case.\medskip 

\n {\bf Question 3.} For vectors $f_1$ and $f_2$ in $L^2_a({\bb B}^m)$, $m>1$, does $[f_1] \simeq [f_2]$ imply anything about the relation of $|f_1|$ and $|f_2|$?

One possible approach would be to try to associate a holomorphic multiplier $\varphi$ with the similarity $X\colon \ [f_1]\to [f_2]$ analogous to the construction in \cite{DM}. That is possible because one can use localization to show that $\dim[f_i]/I_\omega\cdot [f_i]=1$ for $i=1,2$, where $I_\omega$ is the maximal ideal of polynomials in ${\bb C}[\pmb{z}]$ that vanish at $\omega$ in ${\bb B}^m\backslash ({\cl Z}(f_1)\cup {\cl Z}(f_2))$. One might be able to show that $\varphi$ is the quotient of two bounded holomorphic functions on ${\bb B}^m$ using the removeable singularities principle to extend them from ${\bb B}^m\backslash({\cl Z}(f_1) \cup {\cl Z}(f_2))$ to ${\bb B}^m$. However, it is not clear how to connect such a function to $f_1$ and $f_2$.

Although such an argument might seem to show, among other things, that ${\cl Z}(f_1) = {\cl Z}(f_2)$, Richter has pointed out that this relation doesn't hold in general since the multiplier $M_{(z-\lambda)}$ has closed range on $L^2_a({\bb D})$ for $|\lambda|<1$ which implies that $[z-\lambda_1]\simeq [z-\mu]$ for all $|\lambda| <1$ and $|\mu|<1$, $\lambda\ne \mu$. But ${\cl Z}(z-\lambda)= \{\lambda\}$ which shows that the zero sets don't have to be equal. Still the question seems reasonable where the answer might involve the Laplacian of $\log|f_1/f_2|$ in the sense of distributions.

\section{Complemented Submodules}\label{sec6}

\indent

There is another result about the Hardy module on the unit disk whose generalization one can explore for other Hilbert modules. Suppose ${\cl S}$ is a submodule of $H^2({\bb D}) \otimes {\bb C}^2$ that is isomorphic to $H^2({\bb D})$;  one can ask about the quotient module ${\cl Q} = H^2({\bb D}) \otimes {\bb C}^2/{\cl S}$. In particular, is ${\cl Q}$ isomorphic to $H^2({\bb D})$? Simple examples show that it need not be isometrically isomorphic to $H^2({\bb D})$ but it might be similar. 

\begin{thm}\label{thm8}
Let $(\theta_1,\theta_2)$ be a pair of functions in $H^\infty({\bb D})$ such that $|\theta_1(e^{it})|^2 + |\theta_2(e^{it})|^2$ $=1$ a.e. Then ${\cl S} = \{\theta_1f \oplus \theta_2f \in H^2({\bb D})\otimes {\bb C}^2\colon \ f\in H^2({\bb D})\}$ is a submodule such that
\begin{itemize}
\item[\rm (1)] ${\cl Q} \cong H^2({\bb D})$ if and only if $\theta_1$ and $\theta_2$ are constant functions; and
\item[\rm (2)] ${\cl Q}= H^2({\bb D})\otimes {\bb C}^2/{\cl S}$ is similar to $H^2({\bb D})$ if and only if $|\theta_1(z)|^2 + |\theta_2(z)|^2 \ge \vp>0$ for some $\vp>0$ and all $z$ in ${\bb D}$.
\end{itemize}
\end{thm}

\begin{proof}
It is easy to see that the operator $M_z$ on  ${\cl Q}$ is an isometry if and only if ${\cl S}$ is a reducing subspace of $H^2({\bb D})\otimes {\bb C}^2$. This happens only when ${\cl S} = H^2({\bb D}) \otimes {\cl D} \subseteq H^2({\bb D}) \otimes {\bb C}^2$ for some subspace ${\cl D}$ of ${\bb C}^2$.

The result in (2) is  a special case of a result of Sz.-Nagy and Foias \cite{NF}.
\end{proof}

Recently in \cite{DFS}, the authors sought to extend the latter result to the Drury--Arveson space and other related Hilbert modules. However, a full generalization eluded us since we were unable to resolve the following question.\medskip 

\n {\bf Question 4.} Let ${\cl S}$ be a submodule of $H^2_m \otimes {\cl E}$ for some Hilbert space ${\cl E}$ such that
\begin{itemize}
 \item[(1)] ${\cl S} \cong H^2_m\otimes {\cl E}_*$ for some Hilbert space ${\cl E}_*$ and
\item[(2)] there exists a submodule $\widetilde{\cl S}$ of $H^2_m\otimes {\cl E}$ such that
\[
 H^2_m \otimes {\cl E} = {\cl S} \dot{+} \widetilde{\cl S}, \text{ a skew direct}.
\]
\end{itemize}
Does it follow that $\widetilde{\cl S}$ is isomorphic to $H^2_m\otimes {\cl D}$ for some Hilbert space ${\cl D}$?
\medskip

Again, it is easy to show that $\widetilde{\cl S}$ is orthogonal to ${\cl S}$ if and only if ${\cl S} = H^2_m \otimes {\cl E}_*$ for some subspace ${\cl E}_*\subseteq {\cl E}$. An affirmative answer to this question is equivalent to a weakened version of the Beurling--Lax--Halmos Theorem. In particular, one knows, due to McCullough--Trent \cite{MT} and Arveson \cite{Arv}, that $\widetilde{\cl S}$ is the range of a partially isometric multiplier. Unfortunately, it is shown in \cite{DFS} that this map must have a nontrivial null space unless $\widetilde{\cl S}\cong H^2_m\otimes {\cl D}$ for some Hilbert space ${\cl D}$. However, it is possible that $\widetilde{\cl S}$ is the range of a multiplier with closed range and no null space. In that case one can show that $\widetilde{\cl S}$ is similar to $H^2_m\otimes {\cl D}_*$ for some Hilbert space ${\cl D}_*$.

One can ask analogous questions about other quasi-free Hilbert modules but one of the key results used in \cite{DFS} is the lifting theorem which is known to hold only for the Drury--Arveson  space \cite{MT} and closely related Hilbert modules.

Let us conclude with a perhaps surprising result for the one variable case and the related question for the multivariate case.

\begin{thm}\label{thm9}
If ${\cl L}_1$ and ${\cl L}_2$ are submodules of $L^2_a({\bb D})\otimes {\bb C}^2$ so that $L^2_a({\bb D})\otimes {\bb C}^2 = {\cl L}_1\dot{+} {\cl L}_2$, then ${\cl L}_1$ and ${\cl L}_2$ are both isomorphic to $L^2_a({\bb D})$.
\end{thm}

\begin{proof}
The result follows from Lemma 5.2.1 in \cite{JiW} since $M_z$ on $L^2_a({\bb D})$ is in $B_1({\bb D}) \cap(SI)$, where the latter is the set of strongly irreducible operators.
\end{proof}

We don't know if the result holds for the $m>1$ case which we formulate as follows.\medskip

\n {\bf Question 5.} Suppose ${\cl L}_1$ and ${\cl L}_2$ are submodules of $L^2_a({\bb D}^m) \otimes {\bb C}^2$ so that $L^2_a({\bb D}^m) \otimes {\bb C}^2 = {\cl L}_1 \dot{+} {\cl L}_2$. Does it follow that ${\cl L}_1$ and ${\cl L}_2$ are each isomorphic to $L^2_a({\bb D}^m)$?

The problem here is that a priori ${\cl L}_1$ and ${\cl L}_2$ might be very different from $L^2_a({\bb D}^m)$ since that is the case for general submodules of $L^2_a({\bb D}^m)$. The question asks if that is still the case for complemented submodules. This is, of course, the simplest example of a whole family of questions.

\end{document}